\documentclass[13pt,psamsfonts]{amsart}

\usepackage{amsmath}
\usepackage{amsthm}
\usepackage{amssymb}
\usepackage{amscd}
\usepackage{amsfonts}
\usepackage{amsbsy}
\usepackage{enumerate}
\usepackage{graphicx}
\usepackage{calc}
\usepackage{xcolor}
\usepackage{setspace}
\usepackage[dvips]{psfrag}
\usepackage{color}
\usepackage{epsfig}
\usepackage{longtable}

\newtheorem {theorem} {Theorem}

\title[Dynamics of Belousov-Zabotinsky model]
{Dynamics of a $2$-dimensional slow-fast Belousov-Zabotinsky model}

\author[R. Xu,
M. Sun,
X. Zhang]
{Ruihan Xu$^\dag$,
Ming Sun$^\dag$,
Xiang Zhang$^\ddag$ }

\address{$^\dag$School of Mathematical Sciences, Shanghai Jiao Tong University, Shanghai 200240, People's Republic of China}
\address{$^\ddag$School of Mathematical Sciences, MOE--LSC, and CMA-Shanghai, Shanghai Jiao Tong University, Shanghai 200240, People's Republic of China}

\email{xrh-2002@sjtu.edu.cn, 110509863@sjtu.edu.cn, xzhang@sjtu.edu.cn}


\keywords{Belousov-Zhabotinsky differential systems; slow-fast systems; global stability; relaxation oscillation; canard explosion.\\
{\bf Xiang Zhang:} Corresponding author}
\begin{document}
\begin{abstract}
For the reduced two-dimensional Belousov-Zhabotinsky slow-fast differential system, the known results are the existence of one limit cycle and its stability for particular values of the parameters. Here, we characterize all dynamics of this system except one degenerate case. The results include global stability of the positive equilibrium, supercritical and subcritical Hopf bifurcations, the existence of a canard explosion and relaxation oscillation, and the coexistence of one nest of two limit cycles with the outer one originating from the supercritical Hopf bifurcation at one canard point and the inner one from the subcritical Hopf bifurcation at another canard point. This last one is a new dynamical phenomenon.
\end{abstract}

\maketitle

\section{Introduction and statement of the main results}
Chemical oscillatory reactions, particularly the Belousov-Zhabotinsky (BZ) reaction, exemplify the intricate behavior of nonlinear dynamical systems. The discovery by Belousov \cite{BZ_Origins} in 1958, detailing a periodic response within his namesake system, laid the groundwork for this field. Zhabotinsky  \cite{Zhabo1964a, Zhabo1964b} expanded upon Belousov's findings, providing a comprehensive analysis of the kinetics involved in the oxidation of malonic acid.

Despite initial doubts due to its deviation from then-understood thermodynamic principles, the BZ reaction has been thoroughly validated and serves as a crucial model for chemical oscillations. The reaction's mathematical representation, mainly through partial differential equations, captures its spatiotemporal dynamics, as shown by studies including those by Huh et al. \cite{BZ_PDE1} and Toth et al. \cite{BZ_PDE2}, who employed novel methodologies like coevolutionary algorithms for dynamic control and information processing within the BZ reaction.

In 1992, Gy\"{o}rgyi and Field \cite{tv} introduced a three-variable model demonstrating deterministic chaos within the BZ reaction, signifying the fine line between periodicity and chaos. Their model, which adhered closely to experimental observations, was more representative of the BZ reaction's behavior than the simpler Oregonator model that preceded it.

Leonov and Kuznetsov  \cite{HD} in 2013 employed numerical methods to demonstrate limit cycles within a slow-fast BZ system
 \begin{equation}\label{e0}
  \begin{split}
  \epsilon \dot x&=x(1-x)+\frac{f(q-x)}{q+x}y,\\
  \dot y&=x-y,
  \end{split}
\end{equation}
offering new insights into the reaction's dynamics. Barzykina  \cite{CM} in 2020 further elucidated the educational value of the BZ reaction, reinforcing its relevance in a laboratory setting. Llibre and Oliveira \cite{LCB} provided rigorous proof of a unique limit cycle in system \eqref{e0}, thus deepening the theoretical understanding to the dynamics of this system.

This paper extends the examination of system~\eqref{e0}, integrating methods from earlier research ~\cite{SIAP} of explore uncharted aspects of the BZ reaction's dynamics.

\begin{figure}[ht]
  \begin{center}
  \begin{tabular}{ccc}
  \includegraphics[width=5.0cm]{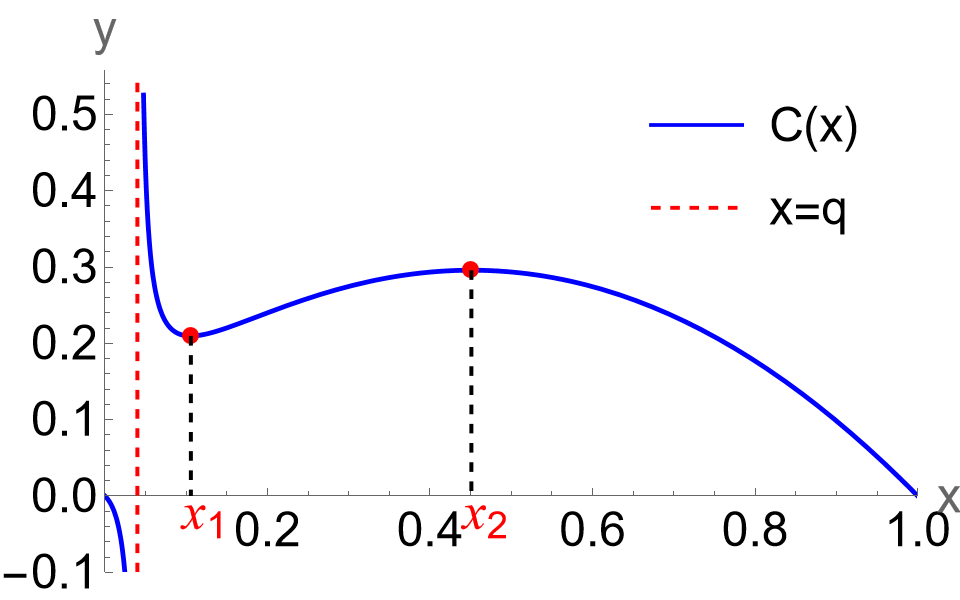} &
  \includegraphics[width=5.0cm]{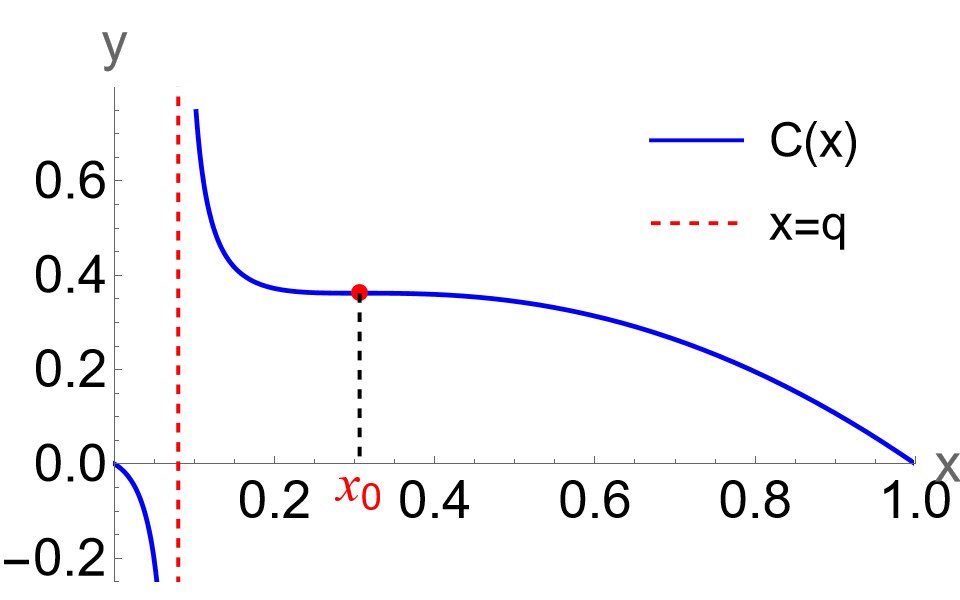} \\
  $(a)$ $0<q<q^*$ & \qquad $(b)$ $q=q^*$ \\
  \includegraphics[width=5.0cm]{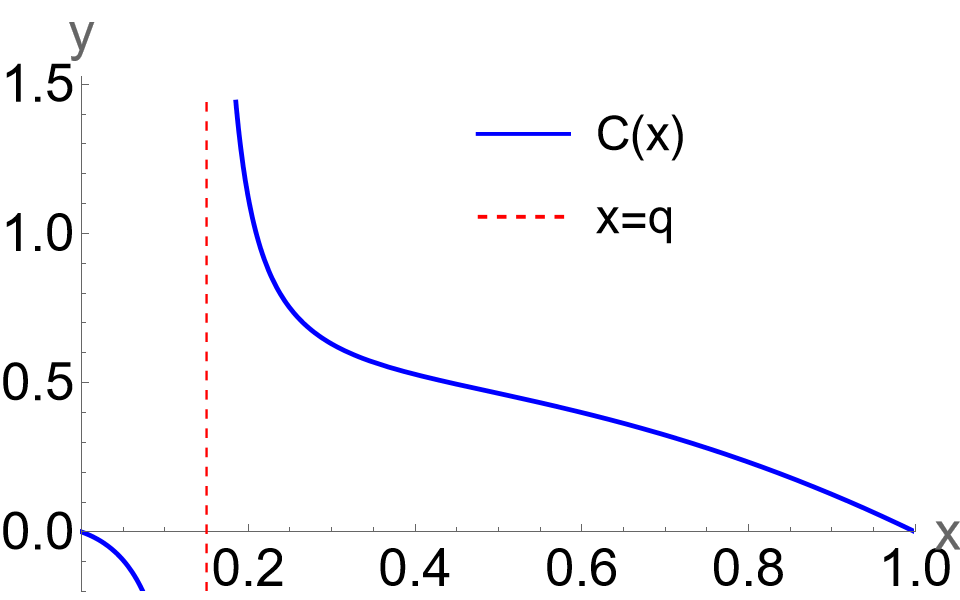} &
  \includegraphics[width=5.0cm]{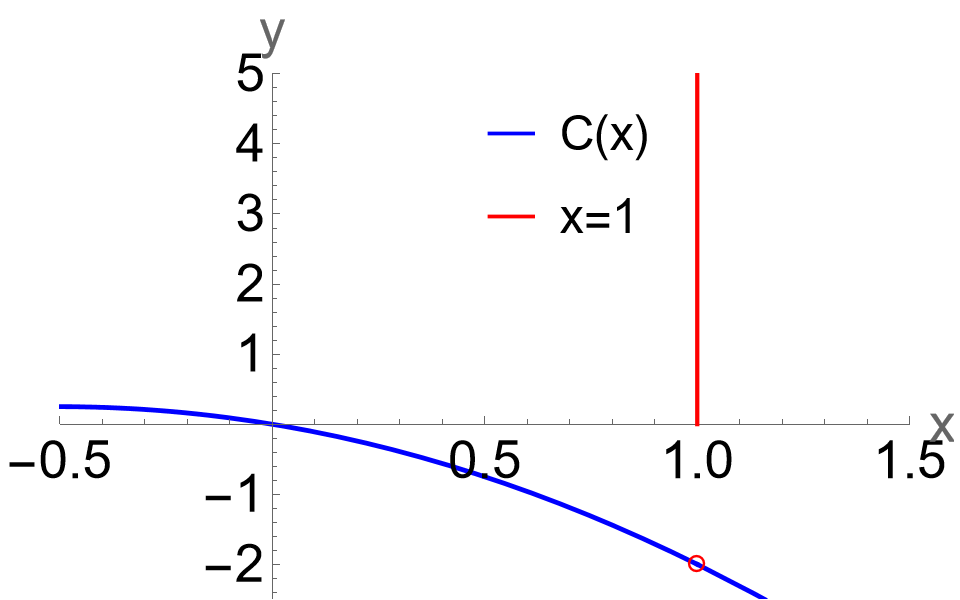} \\
  $(c)$ $q^*<q<1$ &  \qquad $(d)$ $q=1$ \\
  \includegraphics[width=5.0cm]{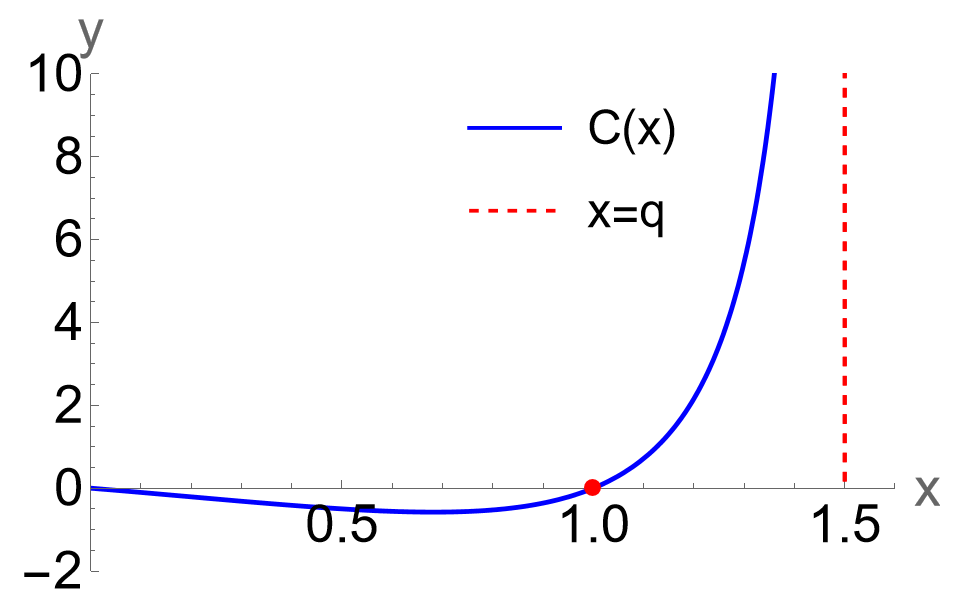} \\
  $(e)$ $q>1$ \\
  \end{tabular}
  \end{center}
  \caption{Shape of the critical curve $S_0$.} \label{dzero}
\end{figure}

In the slow-fast framework, we need to combine the limiting systems of the \textit{slow system} \eqref{e0} and of the \textit{fast system}
\begin{equation}\label{e00}
\begin{split}
\frac{dx}{d\tau}&=x(1-x)+\frac{f(q-x)}{q+x}y,\\
\frac{dy}{d\tau}&=\epsilon(x-y),
\end{split}
\end{equation}
where $\tau=t/\epsilon$ is called \textit{fast time} and $t$ is the \textit{slow time}.
Taking $\epsilon \rightarrow 0$ in the fast system \eqref{e00} and the slow system \eqref{e0}, one gets the \textit{layer system}
\begin{equation}\label{e3}
  \begin{split}
  \frac{dx}{d\tau} &=x\left(1-x \right)+\frac{f(q-x)}{q+x}y,\\
  \frac{dy}{d\tau} &=0,
  \end{split}
\end{equation}
and the \textit{reduced system}
\begin{equation}\label{e4}
  \begin{split}
  0&=x\left(1-x \right)+\frac{f(q-x)}{q+x}y,\\
  \dot y&=x-y.
  \end{split}
\end{equation}
By definition, the set
\[
M_0=\left\{(x,y)\in\mathbb R_+^2| \ x(1-x)+\frac{f(q-x)}{q+x}y=0\right\},
\]
is called the \textit{critical set}, whose geometry deeply affects dynamics of system \eqref{e0}, or equivalently \eqref{e00} when $\epsilon\ne 0$. In the case that $M_0$ is a manifold, it is called a \textit{critical manifold}. The geometric configuration of $M_0$ is the following.
\begin{itemize}
\item[$(a)$] For $q=1$, the critical set $M_0$ is composed of the line $x=1$ and the quadratic curve $y=C_0(x)=-f^{-1}x(x+1)$.
\item[$(b)$] For $q\ne 1$, the critical set $M_0$ has the expression
\[
S_0=\left\{y=C(x):=\frac{x(1-x)(q+x)}{f(x-q)},\ \ x\ne q \right\}.
\]
\end{itemize}
Easy calculations show that
\[
C'(x)=\frac{-2x^3+(1+2q)x^2+2q(q-1)x-q^2}{f(x-q)^2}.
 \]
Set
\[
q^*:=-0.2+0.6\times 2^{\frac{1}{3}}-0.3\times 2^{\frac{2}{3}},
\]
approximately 0.07973, which is a value of $q$ determining the number of critical points of $S_0$. One can check that on the critical curve $S_0$,
$C'(x)$ has no a zero for $q > q^*$, a unique zero saying $x_0$ for $q=q^*$ and that for $0 < q < q^*$, $C'(x)$ has two zeroes, denoted by $x_1$ and $x_2$ with $q<x_1 < x_2<1$, which are marked in Fig.\ref{dzero}$(a)$.

For $0<q<q^*$, $S_0$ has a minimum point $m:=(x_1,C(x_1))$ and a maximum point $M:=(x_2,C(x_2))$, which separate $S_0$ into three subsets:
\begin{align*}
  S_l &= \{(x, y) \in S_0 \mid x < x_1\}, \quad S_m = \{(x, y) \in S_0 \mid x_1 < x < x_2\}, \\
  S_r &= \{(x, y) \in S_0 \mid x > x_1\}.
\end{align*}

Note that the equilibria (if exist) of system \eqref{e0} are always located on the critical curve $S_0$. Their relative locations depend on the values of the parameter $f$. 

The following results characterize the dynamics of system \eqref{e0}.

\begin{theorem}\label{t1}
  For system \eqref{e0} with $f,\ q>0$, and $0<\epsilon\ll 1$, the following statements hold.
  \begin{itemize}
  \item[$(a)$] System \eqref{e0} has exactly one positive equilibrium {in the first quadrant}, saying $E_*=(x_*,x_*)$, which is located either in the strip $q<x<1$ when $q<1$, or in the strip $1<x<q$ when $q>1$, or on the line $x=1$ when $q=1$. 
  \item[$(b)$] If $q\ge q^*$, the unique positive equilibrium is locally stable and so is globally stable in the region $x,y>0$.
  \item[$(c)$] Assume that $0<q<q^*$. Then the critical curve $S_0$ is of $S$-shaped.
  \begin{itemize}
  \item[$(c_1)$] For either $ x_*\in(q, \ x_1]$ or $ x_*\in (x_2,\ 1)$, the equilibrium $E_*$ is locally stable, and so is globally stable .

  \item[$(c_2)$] For any fixed $q\in(0,q^*)$, there exist positive functions $d_1(\epsilon),d_2(\epsilon)=O(\epsilon)$ such that if $x_*\in\left (x_{1\epsilon},\ x_{2\epsilon}\right)$, with $x_{1\epsilon}=x_1+d_1(\epsilon)$ and $x_{2\epsilon}=x_2-d_2(\epsilon)$, the equilibrium $E_*$ is an unstable node or focus. The system has a stable limit cycle, denoted by $\Gamma$, which exhibits relaxation oscillation.

For $x_*= x_{1 \epsilon}$, there is a supercritical Hopf bifurcation at $E_*=E_m:=(x_{1\epsilon},x_{1\epsilon})$.    There exists a $q^{**}\sim 0.05551$ such that at $E_*=E_M:=(x_{2\epsilon},x_{2\epsilon})$ there happens a supercritical Hopf bifurcation when $q>q^{**}$, and a subcritical Hopf bifurcation when $q<q^{**}$. 

  \item[$(c_3)$] Moreover, with monotonic decreasing of $f$, { which leads to monotonic increasing of $x_*$ from $x_1$ to $x_2$}, the limit cycle $\Gamma$ could experience canard explosion as shown below.

For $q^{**}<q<q^*$, $\Gamma$ may undergo a canard explosion, relaxation oscillations and an inverse canard explosion.

For $0<q<q^{**}$, $\Gamma$ may undergo a canard explosion, relaxation oscillations, and a limit cycle of multiplicity two, which is formed via $\Gamma$ coinciding with another limit cycle birthing from the subcritical Hopf bifurcation at $E_*=E_M$.
  \end{itemize}

  \end{itemize}
\end{theorem}

We remark that the known results on system \eqref{e0} are only on the existence of a limit cycle via simulation by Leonov and Kuznetsov \cite{HD} in 2013 and on the existence of at least one limit cycle and its stability for special values of the parameters by Llibre and Oliveira \cite{LCB} in 2022 analytically. Here, we characterize all dynamics of the system except when $q=q^{**}$, which is a degenerate case that we cannot study its dynamics at present. Moreover, the last conclusion in statement $(c_3)$ and the proof to it show that system \eqref{e0} can have one nest of two limit cycles, in which one originates from a supercritical Hopf bifurcation and another one stems from a subcritical Hopf bifurcation. To our knowledge, this is a new dynamical phenomenon not previously found in the published works. 

\begin{figure}[ht]
\begin{center}
\begin{tabular}{c}
\includegraphics[width=10cm]{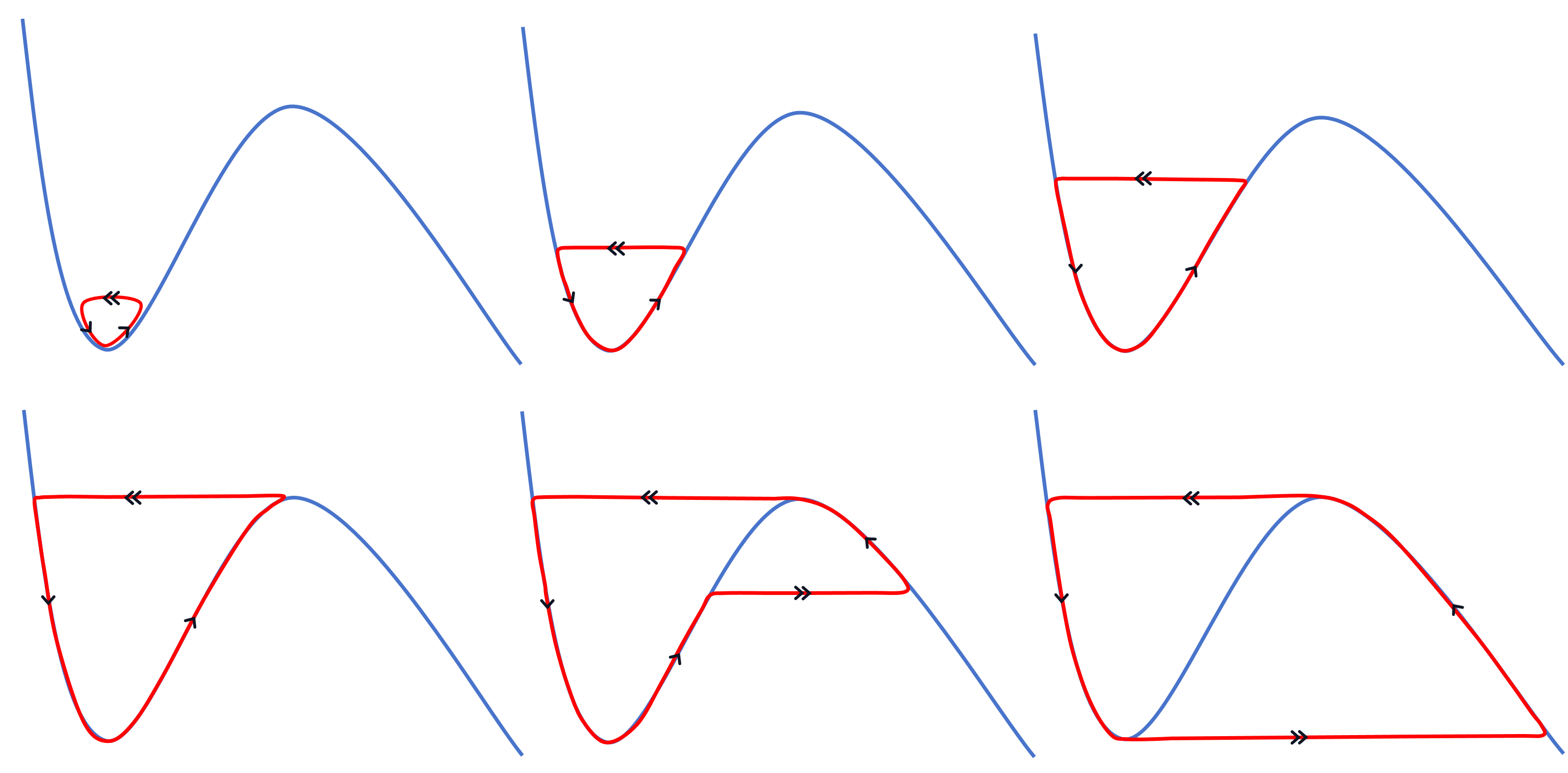}
\end{tabular}
\end{center}
\caption{Illustration of a canard explosion which is a variation of the limit cycle in shape from a small Hopf cycle, a small canard cycle without head, a big canard cycle without head, the transition cycle, a canard cycle with head, and the relaxation oscillation.} \label{fce}
\end{figure}

Recall by definition in \cite{KS2001JDE} that a \textit{canard explosion} is a dynamical phenomenon for slow-fast systems with both a perturbation parameter and a braking parameter. Fixing the perturbation parameter sufficiently small and varying the braking parameter monotonically from a suitable value in an exponentially small interval, the system has a limit cycle birthing from a singular Hopf bifurcation with its shape varying from a small canard cycle to a big canard cycle without head, to a transition cycle, to a canard cycle with head, and finally to a relaxation oscillation. Fig.\ref{fce} illustrates the variation process of these different types of cycles, and Fig. \ref{sim} \ref{sim2} presents numerical simulations of a canard explosion for system \eqref{e00}.

The remaining is the proof of Theorem \ref{t1}.

\section{Proof of Theorem \ref{t1}}

\subsection{Proof of Theorem \ref{t1}$(a)$}

We commence by determining the equilibrium points of system \eqref{e0}. Obviously, any equilibrium is of the form $(x,x)$ with $x$ satisfying the equation
\[
x(1-x) + \frac{f(q - x)}{q + x}x = 0.
\]
It is easy to check that this equation can have a negative root, a zero root, and a positive root with the expression
\[
x_*=\frac{1-q-f+\sqrt{(q+f-1)^2+4q(1+f)}}{2}. 
\]
Only the positive equilibrium $(x_*,x_*)$ is of dynamical interest. Moreover, it is easy to verify that $x_*=1$ if and only if $q=1$; and $x_*<1$ if and only if $q<1$. Further calculations show that if $q<1$ then $x_*>q$; if $q>1$ then $x_*<q$.

In summary, system \eqref{e0}, equivalently \eqref{e00}, has a unique positive equilibrium, denoted by $E_*=(x_*,x_*)$ with $ x_*\in(q, 1)$ for $q<1$, or $x_*\in(1,q)$ for $q>1$, or $x_*=1$ for $q=1$.

\subsection{Proof of Theorem \ref{t1}$(b)$}


We have $q\ge q^*$ by the assumption of the statement. For $q>q^*$ the critical curve $S_0$ strictly monotonically decreases in the interval $(q,1)$. While when $q=q^*$, the critical curve $S_0$ strictly monotonically decreases in the interval $(q,1)$ except at a unique critical point $C_0=(x_0,y_0)$ of the curve.
The Jacobian matrix of system \eqref{e00} at $E_*$ is
\[J_{E_*}=
\begin{pmatrix}
-h(x_*)C'(x_*)& h(x_*)\\
\epsilon & -\epsilon\\
\end{pmatrix}.
\]
where $h(x)=\frac{f(q-x)}{q+x}$. The eigenvalues of the Jacobian matrix are
\begin{equation}\label{eeigen}
\lambda_{1,2}=\frac{-w(x_*,\epsilon) \pm \sqrt{\left(w(x_*,\epsilon)\right)^2-4\left(C'(x_*)-1\right)h(x_*)\epsilon}}{2},
\end{equation}
with $w(x_*,\epsilon)= h(x_*)C'(x_*)+\epsilon $.

For $q^*< q<1$, the positive equilibrium $E_*$ is located in the strip $q<x<1$, and the critical curve $S_0$ monotonically decreases in the region $q<x<1$. This means that $h(x_*)<0$ and $C'(x_*)<0$, and consequently Re$\lambda_{1,2}<0$. So $E_*$ is a stable node or focus.
For $q=q^*$, we also have $h(x_*)<0$ and $C'(x_*)<0$ except $E_*$ is at $C_0$ where $C'(x)=0$. So $E_*$ is also hyperbolic and stable when $x_*\in(q,1)\setminus\{x_0\}$. 

We now adopt the Dulac criterion \cite{Y} to prove the global stability of the equilibrium $E_*$ in the region $\Omega_1:=\{(x,y)|q<x<1,y>0\}$. Set $H(x)=\frac{q+x}{f(q-x)}$, and let $(P(x,y),Q(x,y))=\big(x(1-x)+\frac{f(q-x)}{q+x}y,\epsilon(x-y)\big)$ be the vector field associated to system \eqref{e00}, one gets
$$
\frac{\partial(PH)}{\partial x}+\frac{\partial(QH)}{\partial y}=-C'(x)-\epsilon H(x)>0 \quad q<x<1,y>0,
$$
where we have used the fact that $C'(x)\le 0$ and $H(x)<0$ for $x\in(q,1)$. From the Dulac criterion, one gets that the equilibrium $E_*$ is globally stable in the strip $q < x < 1$, and consequently globally stable in the region $x,y>0$ since $\dot x>0$ on $x\leq q$ and $\dot x<0$ on $x\geq1$.

For $ q=1$, $E_*=(1,1)$. At which system \eqref{e00} has two negative eigenvalues $\lambda_1=-\left(1+\frac{f}{2}\right)$ and $\lambda_2=-\epsilon$. 
{Note that $\dot x = P(x,y)=(1-x)(x+\frac{f}{1+x}y)$. With the fact that $P(x,y)|_{x < 1,y>0}>0$ and $P(x,y)|_{x > 1,y>0}<0$, one gets that all orbits inside the first quadrant will be first attracted to the line $x=1$ and then further attracted to $E_*$. Consequently, $E_*$ is a globally stable node.}


For $ q>1$, the positive equilibrium $E_*$ is located in the strip $1<x<q$, and the critical curve $S_0$ strictly monotonically increases in the region $1<x<q$. This means that $h(x_*)>0$ and $C'(x_*)>0$, and consequently Re$\lambda_{1,2}<0$. So $E_*$ is a stable node or focus. Same as the situation $q^*<q<1$, we have
$$
\frac{\partial(PH)}{\partial x}+\frac{\partial(QH)}{\partial y}=-C'(x)-\epsilon H(x)<0, \quad 1<x<q,\ y>0.
$$
Combined with $\dot x > 0$ on $x\leq1$ and $\dot x < 0$ on $x\geq q$, it follows that $E_*$ is globally stable in the region $x,y>0$. Statement $(b)$ follows.


\subsection{Proof of Theorem \ref{t1}$(c)$}

Now we have $0<q<q^*$ by the assumption of the statement. The critical curve $S_0$ is that given in Fig.\ref{dzero}$(a)$.

\noindent $(c_1)$ Since $x_*\in (q, x_1]$ or $x_*\in[x_2,1)$, one has $h(x_*)<0$ and $C'(x_*)<0$ for $x_*\ne\{x_1,x_2\}$ and $C'(x_*)=0$ for $x_*\in\{x_1,x_2\}$. It follows from \eqref{eeigen} that $E_*$ is hyperbolic and stable for $x_*\in (q, x_1)$ or $x_*\in(x_2,1)$. 

{We first discuss the global stability of $E_*$ with $x_*\in (q,x_1)$, i.e. $E_*\in S_l$.
Note that the two extreme points $m$ and $M$ of the critical curve $S_0$ are both jump points \cite{KSSIAM2001,KS2001JDE}, and that $S_0\setminus\{m,M\}$ is normally hyperbolic.
For any positive orbit $O_\Gamma$ in the region $x>q, \ y>0$ whose initial point $I_\Gamma$ is not on $S_0$, we have the next conclusions:
 \begin{itemize}
\item If $I_\Gamma$ is above the critical curve $S_0$, 
 it follows from the Fenichel invariant manifold theory \cite{IM} that $O_\Gamma$ will move along one of the fast orbit of the layer system to an $O(\varepsilon)$ neighborhood of either $S_l$ or of $\{M\}\cup S_r$. In the first case, $O_\Gamma$ will follow $S_l$ moving to $E_*$.
 In the second case, by the jump point theory \cite{KS2001JDE} $O_\Gamma$ will jump from $M$ to $S_l$ and then follow $S_l$ to $E_*$. In the last case, $O_\Gamma$ will follow $S_r$ moving up to a neighborhood of $M$, then jump from $M$ to $S_l$ and finally be attracted to $E_*$.
 \item If $I_\Gamma$ is under the critical curve $S_0$,
 again by the Fenichel invariant manifold theory \cite{IM} and the jump point theory \cite{KS2001JDE} $O_\Gamma$ will be first attracted to an $O(\varepsilon)$ neighborhood of either $S_l\cup\{m\}$ or of $S_r$. The same arguments as used in the last case show that $O_\Gamma$ will be finally attracted to $E_*$.
\end{itemize}
This proves that $E_*$ is globally stable. 

For $x_*\in(x_2,1)$, that $E_*$ is globally stable can be verified in a similar way as above. The proof of statement $(c_1)$ is done. 
}

\noindent $(c_2)$ For $x_*\in (x_1,x_2)$, one has $C'(x_*)> 0$. Sine $q<q^*<1$, it follows $x_*>q$, and so $h(x_*)<0$. Since $\epsilon>0$ is sufficiently small, by the monotonicity of $w(x,\epsilon)$ in $x$ near $x=x_1$ and $x=x_2$ it follows that there exist unique functions $d_1(\epsilon)$ and $d_2(\epsilon) $ which are both positive when $0<\epsilon\ll 1$ and both vanish when $\epsilon=0$, for which the following hold.
\begin{itemize}
\item $w(x_*,\epsilon)=0$ when $x_*=x_1+d_1(\epsilon)$ and $x_*=x_2-d_2(\epsilon)$.

\item $w(x_*,\epsilon)>0$ when $ x_*\in (x_1,\ x_1+d_1(\epsilon))$ and $x_*\in( x_2-d_2(\epsilon),x_2)$.

\item $w(x_*,\epsilon)<0$ when $ x_*\in ( x_1+d_1(\epsilon), x_2-d_2(\epsilon))$.

\end{itemize}
The equation that $x_*=x_1+d_1(\epsilon)$ satisfies gives
\[
(6x^2-2(1+2q)^2x-2q(q-1))d_1+O(d_1^2)+((x+d_1)^2-q^2)\epsilon=0,
\]
where we have used the fact that $C'(x_1)=0$. This implies that $d_1(\epsilon)=O(\epsilon)$.

\noindent \textit{Case }1. {For $ x_*\in [x_1,\ x_1+d_1(\epsilon))$ or $x_*\in( x_2-d_2(\epsilon),x_2]$, since $0<\varepsilon\ll 1$ one has $-w(x_*,\epsilon)<0$. The equilibrium $E_*$ is a stable focus.} 
Similar to the proofs as those in $(c_1)$, by the Fenichel theory \cite{IM} and the jump point theory \cite{KS2001JDE} all orbits starting from the region $x>q,\ y>0$ will be attracted to an $O(\varepsilon)$ neighborhood of either $m$ when $ x_*\in [x_1,\ x_1+d_1(\epsilon))$ or $M$ when 
$x_*\in( x_2-d_2(\epsilon),x_2]$. Then by \cite[Theorem 3.3$(i)$]{KS2001JDE} the orbits will be further attracted to $E_*$, because the neighborhood of the canard point in Theorem 3.3$(i)$ of \cite{KS2001JDE} is $O(\varepsilon^\nu)$ with $\nu\in (0,1)$ and it is larger than $O(\varepsilon)$. So, $E_*$ is also a global attractor in the dynamically interesting region.

\noindent \textit{Case }2. For $ x_*\in ( x_1+d_1(\epsilon), x_2-d_2(\epsilon))$, one has $-w(x_*,\epsilon)>0$. { 
Note that the eigenvalues $\lambda_{1,2} $ of $E_*$ as an equilibrium of system \eqref{e00} 
satisfy
$$\lambda_1\cdot \lambda_2 = -\epsilon(1-2x_*+h'(x_*)y_*-h(x_*))>0,
$$
and $\textrm{Re}\lambda_{1,2}=-w(x_*,\epsilon)>0$, 
where we have used $h'(x_*)=-\frac{2fq}{(q+x)^2}<0$ and $h(x_*)=1-x_*$. So, the equilibrium $E_*$ is an unstable focus or node.} 

Since $q<1$ and all orbits in the first quadrant will positively get into the region $q<x<1$ and $y>0$, one can easily construct an outer boundary satisfying the Poincar\'e-Bendixson annular theorem \cite{Y}, which verifies existence of a stable limit cycle surrounding $E_*$.

\noindent \textit{Case }3. For $x_*=x_1+d_1(\epsilon)$ and $x_*=x_2-d_2(\epsilon)$, one has $w(x_*,\epsilon)=0$. The equilibrium $E_*$ is a weak focus. Instead of computing the focal quantity of system \eqref{e00} at $E_*$, we adopt the singular perturbation approach via Theorem 3.3 of \cite{KS2001JDE} by Krupa and Szmolyan. To this aim we need to verify their conditions \((A1)\) through \((A4)\), and calculate the quantity \(A\) defined therein.

Through the preceding analysis, we have verified the following facts:
The critical curve \(S_0\) is of \(S\)-shaped in the half plane $x>q$ with the minimum point at $m=(x_1, y_1)$ and the maximum point at $M=(x_2, y_2)$, where \(q < x_1 < x_2 < 1\), which separate $S_0$ in three parts \(S_l, S_m, S_r\). The minimum and maximum points are both fold points because
\[
C'(x)<0\quad \mbox{ on }\  S_{l}\cup S_{r}, \qquad \mbox{ and } \qquad C'(x)>0 \quad \mbox{ on } \ S_{m},
\]
and moreover
\begin{align*}
\frac{d C}{d x}(x_{1})&= 0, \quad \frac{d^2 C}{d x^2}(x_{1})>0,\quad \frac{d C}{d x}(x_{2})= 0, \quad \frac{d^2 C}{d x^2}(x_{2})<0,\\
g(x_*,x_*)&=0,\quad \frac{\partial g}{\partial x}(x_*,x_*)\neq 0, \quad \mbox{when } \ \ x_*=x_{1},\ x_2,\\
\end{align*}
where $g(x,y)=x-y$. 
Additionally, writing the reduced system in the form
\[
C'(x)\dot x=x-C(x)
\]
gives the following facts:
\begin{itemize}
\item When $x_*=x_1$, $\dot{x}>0$ on $S_{l}\cup \{m\}\cup S_m$ and $\dot{x}<0$ on $S_{r}$.
\item When $x_*=x_2$, $\dot{x}>0$ on $S_{l}$ and $\dot{x}<0$ on $ S_m \cup \{M\}\cup S_{r}$.
\end{itemize}
These arrive at the next conclusions.
\begin{itemize}
\item If $x_*=x_1$, then $m:=(x_1,y_1)$ is a canard point, and $M:=(x_2,y_2)$ is a jump point.
\item If $x_*=x_2$, then $m $ is a jump point, and $M$ is a canard point.
\end{itemize}
Hence, the assumptions $(A1)$--$(A4)$ in \cite{KS2001JDE} hold.

Finally, we calculate the quantity $A$ defined in \cite{KS2001JDE} and determine its sign. According to \cite{LCB}, after the time rescaling $d\tau=(x+q)ds$ system \eqref{e00} can be written in an equivalent way as the next polynomial differential system
\begin{equation}\label{e1}
  \begin{split}
  x'&=x\left(1-x \right)(q+x)+f(q-x)y,\\
  y'&=\epsilon(x-y)(q+x),
  \end{split}
\end{equation}
where the prime represents the derivative with respect to the new time $s$. Further replacing $fy$ by $y$ gives
\begin{equation}\label{e6}
  \begin{split}
  x'&=x\left(1-x \right)(q+x)+(q-x)y,\\
  y'&=\epsilon(fx-y)(q+x),
  \end{split}
\end{equation}
Translating $(x_*,y_*)=(x_i,y_i)$, $i=1,2$, to the origin $(0,0)$, and writing the resulting system in the normal form system provided in \cite{KSSIAM2001}, we get from system \eqref{e6} the next system
\begin{equation}\label{enorm1}
\begin{split}
  \dot u&=-v\tilde{h}_{1}(u, v,\lambda,\epsilon)+u^2\tilde{h}_{2}(u,v,\lambda,\epsilon)+\epsilon \tilde{h}_{3}(u,v,\lambda,\epsilon),\\
  \dot v&=\epsilon (u\tilde{h}_{4}(u,v,\lambda,\epsilon)-\lambda \tilde{h}_{5}(u,v,\lambda,\epsilon)+v \tilde{h}_{6}(u,v,\lambda,\epsilon)),
\end{split}
\end{equation}
where $u = x-x_*$, $v=y-y_*$, $\lambda = f_*-f $, $f_*x_*=y_*$, and
\begin{align*}
  \tilde{h}_{1}(u,v,\lambda,\mu)&=x_*-q+u,\\
  \tilde{h}_{2}(u,v,\lambda,\mu)&=\frac{1}{2}C''(x_*)(x_*-q)
  +\left(\frac{1}{2}C''(x^*)+\frac{1}{6}C'''(x^*)(x^*-q)\right)u+O(u^2),\\
  \tilde{h}_{3}(u,v,\lambda,\mu)&=0,\\
  \tilde{h}_{4}(u,v,\lambda,\mu)&=f_*(q+ x_*)+f_*u,\\
  \tilde{h}_{5}(u,v,\lambda,\mu)&=qx_*+x^2_*+(q+2x_*)u+O(u^2),\\
  \tilde{h}_{6}(u,v,\lambda,\mu)&=-q-x_*-u,
\end{align*}
where $ f_* = C(x_*)/x_*$.

As in (3.3) of \cite{KSSIAM2001} in order for the constant terms of $\tilde{h}_{1},\tilde{h}_{2},\tilde{h}_{4},\tilde{h}_{5}$ to be equal to 1, we take the rescaling
$$
u\rightarrow \alpha u,v \rightarrow \beta v, \lambda \rightarrow \eta \lambda, t \rightarrow \xi t
$$
with
\[
\begin{array}{ll}
\alpha = \dfrac{2 \sqrt{(x_*^2-q^2)f_* }}{C''(x_*)(x_*-q)}, & \beta = \dfrac{2 f_*(q+x_*)}{ C''(x_*)(x_*-q)}, \\
\eta = \dfrac{ 2 f_*\sqrt{(x_*^2-q^2) f_*}}{ C''(x_*)(x_* -q ) x_*}, & \xi = \dfrac{1}{\sqrt{(x_*^2-q^2) f_*}}.
\end{array}
\]
Then we get system \eqref{enorm1} with $\tilde h_j$'s replaced by the next $h_j$'s
\begin{align*}
  h_1 &= 1+\frac{1}{x_*-q}\alpha u, \\
  h_2 &= 1+\frac{C''(x_*)+\frac{1}{3}C'''(x_*)(x_*-q)}{C''(x_*)(x_*-q)}\alpha u+ O(u^2),\\
  h_4 &= 1 + \frac{1}{q+x_*}\alpha u,\\
  h_6 &= (-q-x_*-\alpha u) \xi.
\end{align*}
Here, we do not specify $h_5$.
According to the calculations obtained in \cite{CZ2021}, one has
\begin{align*}
a_2&=\frac{\partial h_1}{\partial u}|_{\mathcal K}=\frac{\alpha}{x_*-q},\\
a_3&=\frac{\partial h_2}{\partial u}|_{\mathcal K}=\left(\frac{1}{x_*-q}+\frac{1}{3}\frac{C'''(x_*)}{C''(x_*)}\right)\alpha,\\
a_4&=\frac{\partial h_4}{\partial u}|_{\mathcal K}=\frac{\alpha}{q+x_*},\qquad
a_5= h_6|_{\mathcal K}= -( q+x_*)\xi.
\end{align*}
where $\mathcal K=(u,v,\lambda,\epsilon)=(0,0,0,0)$. Then adopting the notations in \cite{KSSIAM2001} gives
\begin{align*}
A=-a_2+3a_3-2a_4-2a_5.
\end{align*}
So one gets
\[
  A(x_*) = \left(\frac{2}{x_*-q}+\frac{C'''(x_*)}{C''(x_*)}-\frac{2}{x_*+q}+\frac{x_* C''(x_*)}{C(x_*)}\right)\frac{2 \sqrt{(x_*^2-q^2) f_*}}{C''(x_*)(x_*-q)}.
\]
Recall that $f_*=C(x_*)/x_*$.

\begin{figure}[ht]
  \begin{center}
  \begin{tabular}{ccc}
  \includegraphics[width=5.5cm]{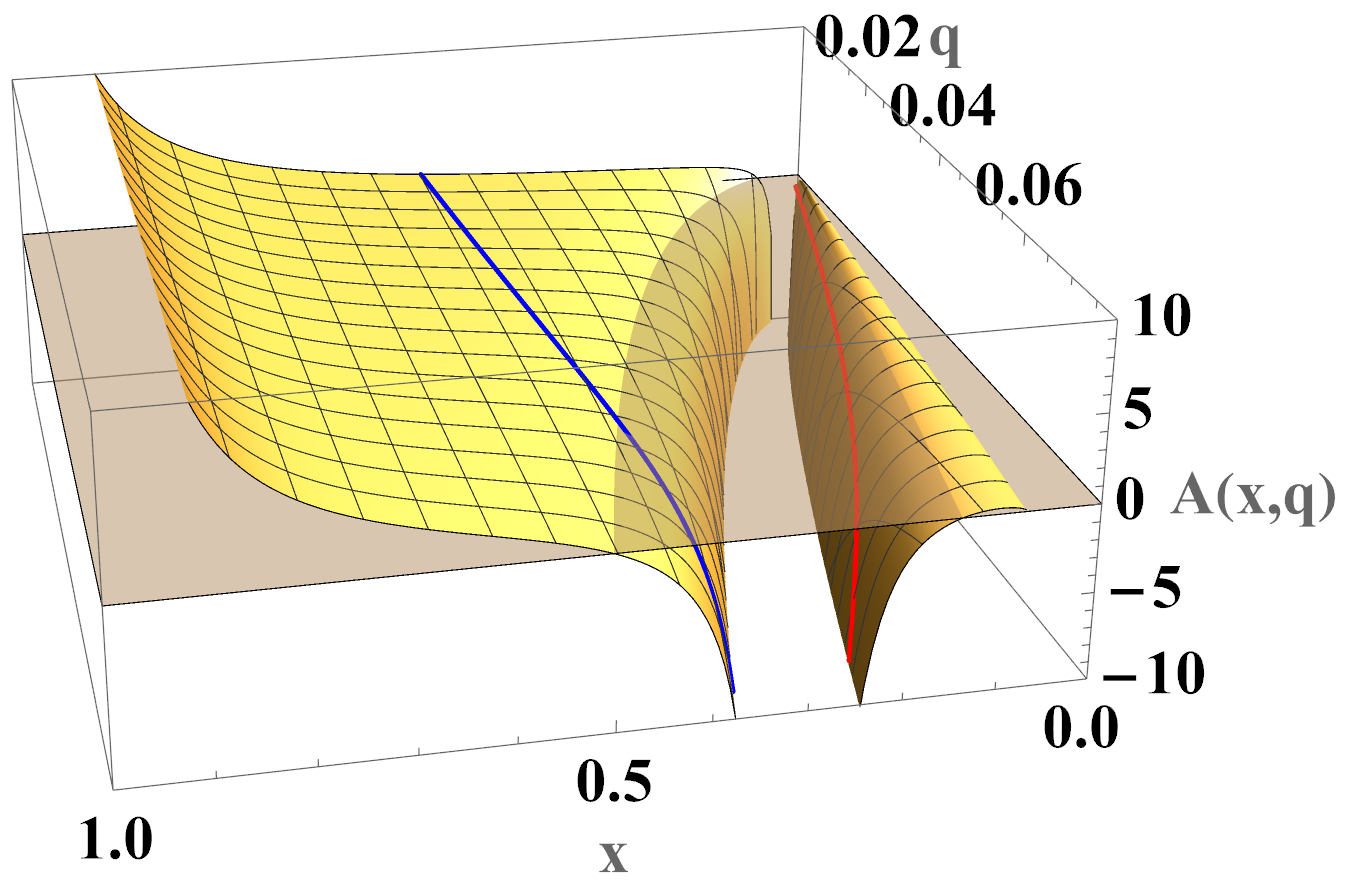} & \includegraphics[width=5.5cm]{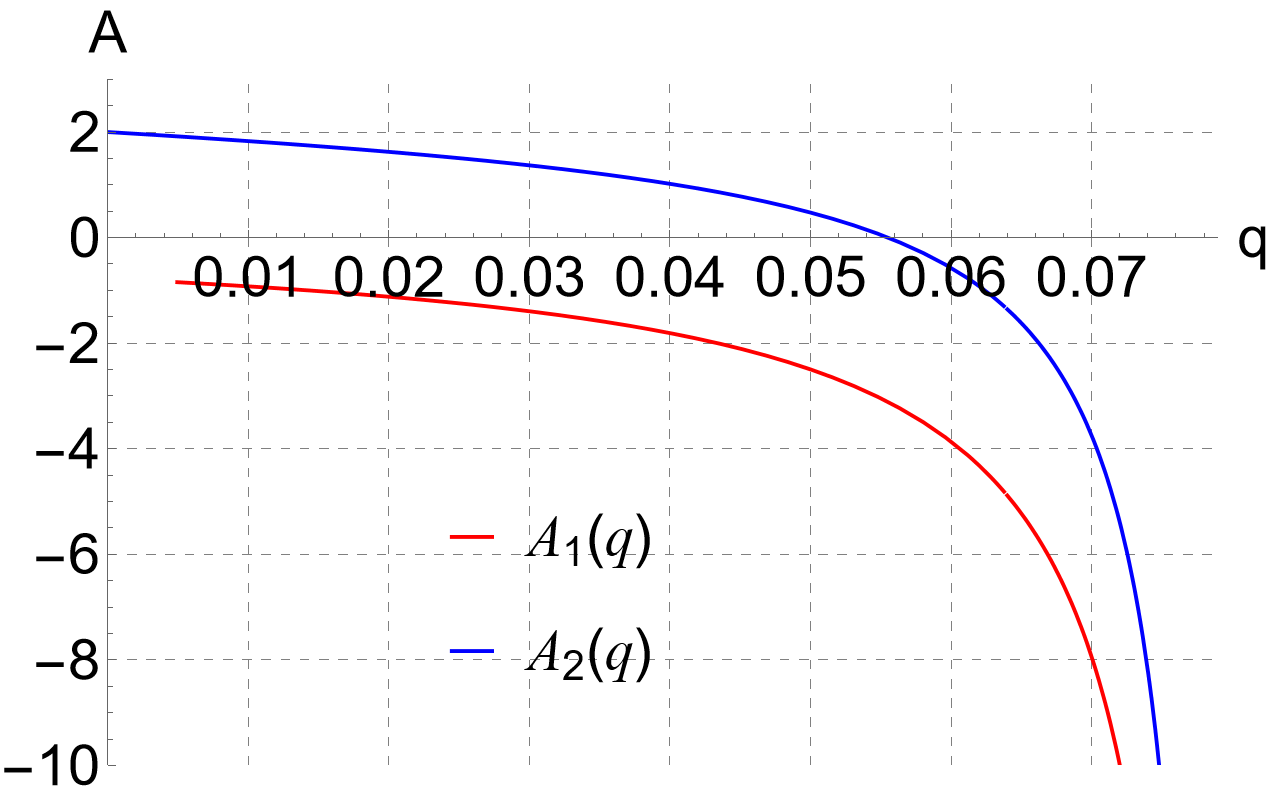}\\
  $(a)$ plot of $A(x,q)$ & $(b)$ curve of $A(q)$
  \end{tabular}
  \end{center}
  \caption{$(a)$ The graph of $ A$ as a function of $(x,q)\in (0,1)\times(0,q^*)$. The right red curve is the graph of $A$ when $x=x_*=x_1$ with $q$ varying from $0$ to $q^{**}$, and the left blue curve is the graph of $A$ when $x=x_*=x_2$ with a variation of $q$. The right one always means $A<0$, and the left one implies that $A$ changes its sign once at a value $q=q^{**}\sim 0.05551$. $(b)$ The red line is the curve of $A_1(q):=A(x_1(q),q)$ and the blue line is the curve of $A_2(q):=A(x_2(q),q)$, since $x_1$ and $x_2$ are uniquely decided by $q$.}
  \label{A}
\end{figure}

The sign of $A$ is hard to determine analytically. Instead, we numerically draw the graph of $A$ in $(x,q)$ in the three-dimensional space. Fig.\ref{A} $(a)$ simulates the graph with $(x,q)\in (0,1)\times(0,q^*)$, where the left blue curve represents the values of $A$ at $x_*=x_2 $, and the right red one denotes the values of $A$ at $x_*=x_1$ with $q$ varying from $0$ to $q^*$. Recall that $x_1$ and $x_2$ are both functions of $q$. Fig.\ref{A} $(b)$ further exhibits the projection of the two curves on the $(A,q)$ plane.

According to Fig. \ref{A}, the value of $A$ at the minimum point $(x_*,x_*)=m=(x_1,x_1)$, denoted by $A_m$, is always negative. Whereas the value of $A$ at the maximum point $(x_*,x_*)=M=(x_2,x_2)$, denoted by $A_M$, changes its sign at $q=q^{**}\sim 0.05551$ from negative to positive when $q$ decreases from $q^*$ to $0$.
Summarizing the above verification gives
\begin{itemize}
\item at $E_*=m$, one has $A=A_m<0$
\item at $E_*=M$, one has $A=A_M<0$ for $q\in (q^{**}, q^*)$ and $A=A_M>0$ for $q\in (0, q^{**})$.
\end{itemize}
According to the asymptotic expression of the first Lyapunov quantity of system \eqref{e6} at $E_*$ near $m$ or $M$, the extreme point, one has the facts: $A<0$ implies that the first Lyapunov quantity of system \eqref{e6} at $E_*$ is negative; $A>0$ means that the first Lyapunov quantity of system \eqref{e6} at $E_*$ is positive. In case $A=0$, i.e. $q=q^{**}$, the situation will be much more complicated, and it will be postponed for future study.

This proves statement $(c_2)$.

\noindent $(c_3)$. {Set $f_1=C(x_1)/x_1$ and $f_2=C(x_2)/x_2$. Then $f_1>f_2$. With the help of system \eqref{e6}, one gets that when $f=f_1$, $m$ is a canard point, and $M$ is a jump point. Whereas when $f=f_2$, $m$ is a jump point, and $M$ is a canard point. When $f$ decreases from $f_1$, the equilibrium $E_*$ will move rightward along $S_m$ from $m$. By Theorems 3.1, 3.2 and 3.3 of \cite{KS2001JDE}, for fixed $q$ and $0<\epsilon\ll 1$ there exist functions $f_{hm}(\epsilon)$ and $f_{cm}(\epsilon)$ satisfying $f_{hm}(\epsilon)>f_{cm}(\epsilon)$ such that
\begin{itemize}
\item[$(R_1)$] $E_m$ is hyperbolic and stable when $f>f_{hm}(\epsilon)$, and it experiences a supercritical Hopf bifurcation when $f$ decreases and passes $f_{hm}(\epsilon)$. Then, a stable limit cycle births via the Hopf bifurcation. Denoted by $\Gamma_{1h}$
\item[$(R_2)$] The limit cycle will undergo a canard explosion for $|f-f_{cm}(\epsilon)|<e^{-\nu_1/\epsilon}$ with $\nu_1>0$. That is, the Hopf cycle expands to small canard cycles without a head, to big canard cycles without a head, to a transition cycle, to canard cycles with a head, and finally to relaxation oscillation, denoted by $\Gamma_{1r}$.
\item[$(R_3)$] In case $q > q^{**}$, since $A=A_M<0$, there exist functions $f_{hM}(\epsilon)$ and $f_{cM}(\epsilon)$ satisfying $f_{hM}(\epsilon)<f_{cM}(\epsilon)$ such that when $f$ decreases such that $|f-f_{cM}(\epsilon)|<e^{-\nu_2/\epsilon}$ with $\nu_2>0$, the relaxation oscillation cycle $\Gamma_{1r}$ will experience an inverse canard explosion and finally shrink to $E_*$ when $f=f_{hM}(\epsilon)$. Then, the system will have no limit cycle for $f<f_{hM}(\epsilon)$.
\item[$(R_4)$] In case $q< q^{**}$, since $A=A_M>0$, when $f$ decreases from $f_{cm}$ to a neighborhood of $f_2$, the relaxation oscillation cycle $\Gamma_{1r}$ via the canard explosion in $(R_2)$ cannot undergo an inverse canard explosion and shrink to the equilibrium $E_m$ near $M$. Instead by Theorems 3.1, 3.2 and 3.4 of \cite{KS2001JDE} there exist functions $f_{hMd}(\epsilon)$ and $f_{cMd}(\epsilon)$ satisfying $f_{hMd}(\epsilon)>f_{cMd}(\epsilon)$ such that when $f$ continuously decreases and is equal to $f_{hMd}$ there happens a subcritical Hopf bifurcation at $E_*$, which gives birth to an unstable limit cycle, saying $\Gamma_{2h}$. With continuous decreasing of $f$ such that $|f-f_{cMd}(\epsilon)|<e^{-\nu_3/\epsilon}$ with $\nu_3>0$, the unstable Hopf cycle $\Gamma_{2h}$ rapidly expands and coincides with the canard cycle deformed from $\Gamma_{1r}$. It becomes a limit cycle of multiplicity $2$, and then system \eqref{e0} has no limit cycle. Here, the multiple $2$ limit cycle bifurcation happens.
\end{itemize}
}

It completes the proof of Theorem \ref{t1}. \qed

\subsection{Simulating the dynamical phenomena}

This section presents numerical simulations via Python on the dynamical phenomena obtained above by our theoretical analysis, i.e., the existence of canard explosion, relaxation oscillation, and multiple $2$ limit cycle, focusing on the situation where $0<q<q^*$ and the critical curve $S_0$ has two critical points at $x=x_{1}$ and $x=x_2$. Recall that the variation of the parameter $q$ leads to a change of the sign of the quantity $A$ when the unique positive equilibrium is located at the maximum point of the critical curve $S_0$, which contributes to different dynamics caused by this sign variation.


First, we take $\epsilon=0.0001, q=0.07\in (q^{**},q^*)$, where the signs of the quantity $A$ at two critical points of $S_0$ are -19.69 and -7.08, both negative. We set the initial point at $(0.4,\frac{0.3379}{f})$ for each $f$. Slightly changing $f$ from 1.536 to 0.89, we get the whole Fig.\ref{sim}, where the solid {green curve} represents the forward orbit from the initial point while the dotted {orange curve} denotes the backward orbit from the same initial point. The left column represents the variation in the limit cycle which undergoes the supercritical singular Hopf bifurcation at the local minimum point $m$ of $S_0$, the canard explosion, and the relaxation oscillation, while the right column illustrates the alteration in the limit cycle from the relaxation oscillation to canard cycles with head, to canard cycles without head, and finally shrinks to the singular Hopf bifurcation point at the local maximum point $M$ of $S_0$.

\begin{longtable}{cc}
  \caption{Numerical simulation in the case that the signs of $A$ are both negative at $E_*=m$ and $E_*=M$} \label{sim} \\
  \endfirsthead \includegraphics[width=0.45\textwidth]{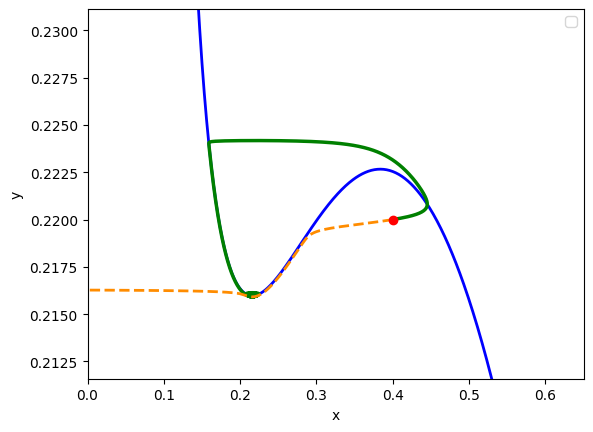} & \includegraphics[width=0.45\textwidth]{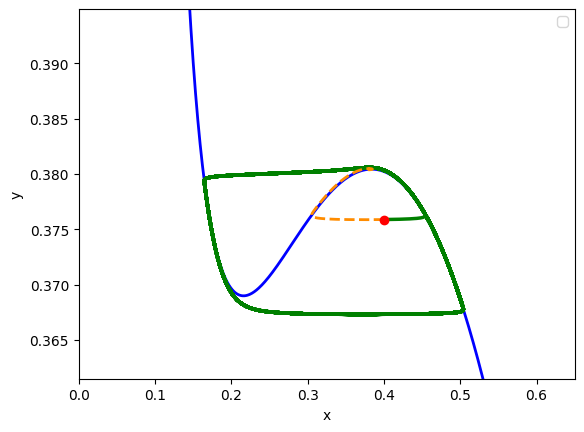} \\
  $(a_1)$\ \  $f=1.536$ & $(b_1)$\ \ $f=0.899$ \\
  \includegraphics[width=0.455\textwidth]{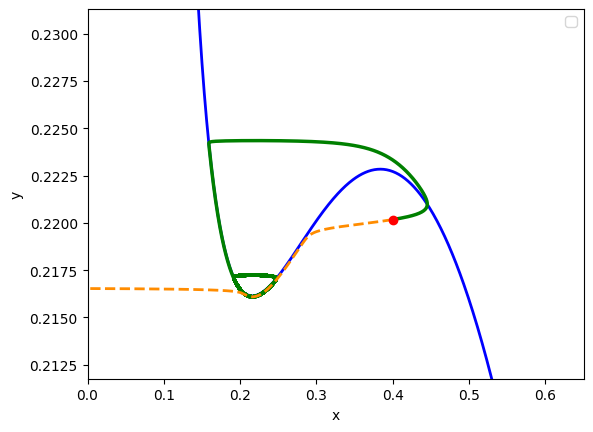} & \includegraphics[width=0.45\textwidth]{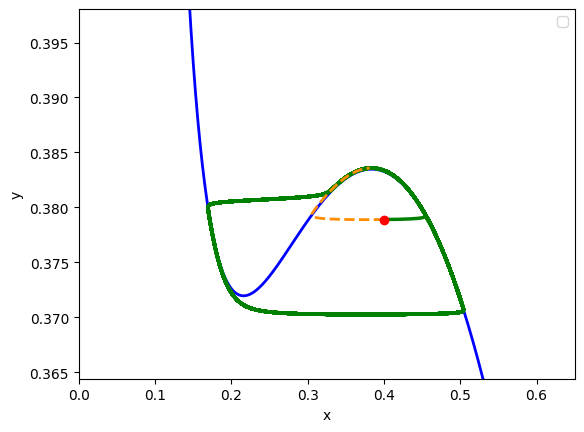} \\
  $(a_2)$\ \  $ f=1.53475666$ & $(b_2)$\ \  $f=0.8918734654$ \\
  \includegraphics[width=0.45\textwidth]{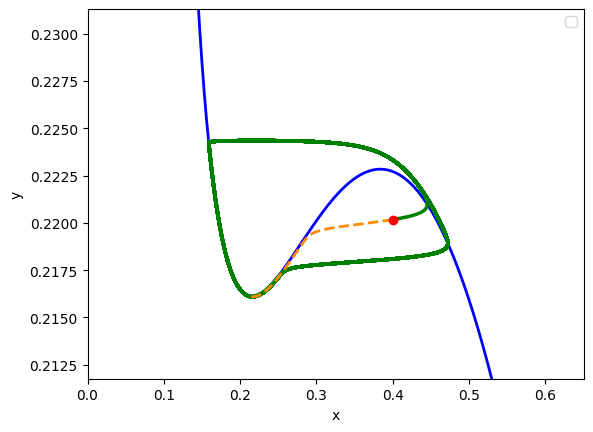} & \includegraphics[width=0.45\textwidth]{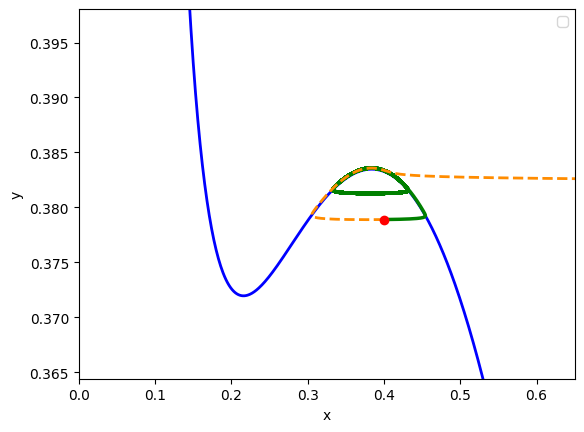} \\
  $(a_3)$\ \  $ f=1.534756651273$ & $(b_3)$\ \  $f=0.891873465$ \\
  \includegraphics[width=0.45\textwidth]{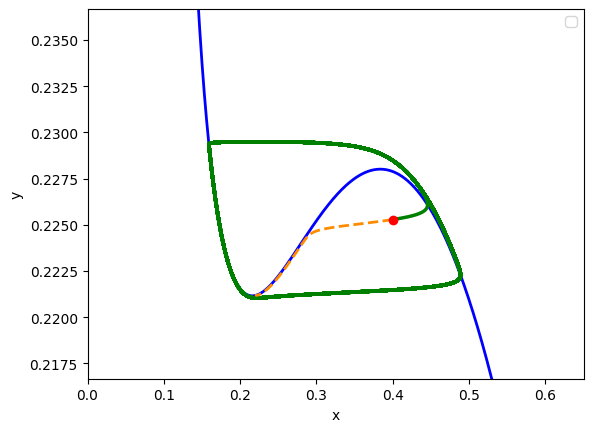} & \includegraphics[width=0.45\textwidth]{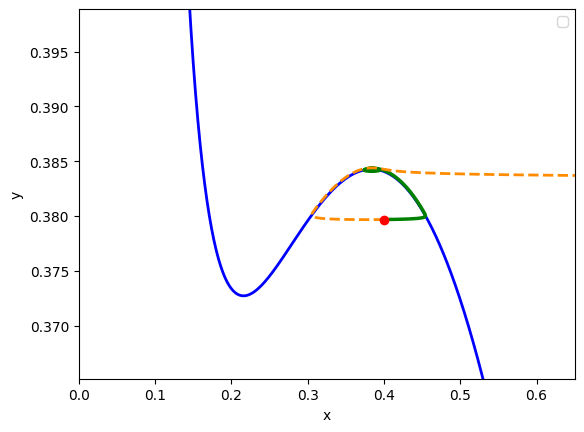} \\
  $(a_4)$\ \  $ f=1.5$ & $(b_4)$\ \  $f=0.89$ \\
\end{longtable}

Next we take $\epsilon=0.01$, $q=0.02\in(0,q^{**})$, where the sign of the quantity $A$ is -3.32 at the local minimum point $m$ but is 2.45 at the local maximum point $M$. { We set the initial point at $(0.4,\frac{0.35}{f})$ for $f=2.18$ and $(0.4,\frac{0.35}{f})$ for other $f$ to illustrate the orbit better}. Slightly changing $f$ (from 2.18 to 0.5), we get the dynamics shown in Fig.\ref{sim2}, which are described as follows.
\begin{itemize}
\item From $(c_1)$ to $(c_4)$ of Fig.\ref{sim2} exhibits the canard explosion starting from the canard point at $E_*=m$.
\item In $(c_5)$ of Fig.\ref{sim2} the unstable equilibrium $E_*$ moves along $S_m$ to the neighborhood of the maximum point $M$ of $S_0$ and at this time, the positive orbit (green one) has the relaxation oscillation as its limit, while the negative orbit (orange one) has the equilibrium $E_*$ as its limit.
\item In $(c_6)$ of Fig.\ref{sim2}, the subcritical singular Hopf bifurcation happens, which gives birth to an unstable limit cycle. The negative orbit has this limit cycle as its limit.
\item From $(c_6)$ to $(c_8)$ of Fig.\ref{sim2} exhibits the shrinking of the outer stable limit cycle and the expansion of the inner unstable limit cycle. With $f$ further decreasing from the value in $(c_8)$, the two limit cycles should coincide and become a double one. But this phenomenon is too sensitive with the value of $f$ to simulate it. 
\item In $(c_9)$ and $(c_{10})$, the system has no limit cycle, and the equilibrium $E_*$ is hyperbolic and stable. The positive orbit has its limit at $E_*$, and the negative orbit passes $M$ and goes to infinity at the positive $x$-direction.
\end{itemize}

We remark that here we take $\epsilon=0.01$ but not $0.0001$, because such a small $\epsilon>0$ will
cause the orbit to be too sensitive to simulate the two limit cycles.

\begin{longtable}{cc}
  \endfirsthead \includegraphics[width=0.45\textwidth]{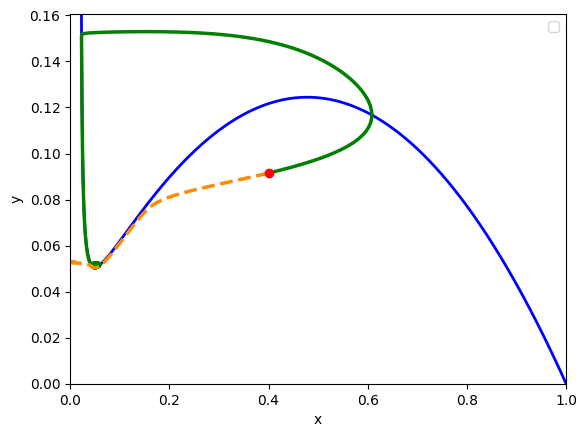} & \includegraphics[width=0.45\textwidth]{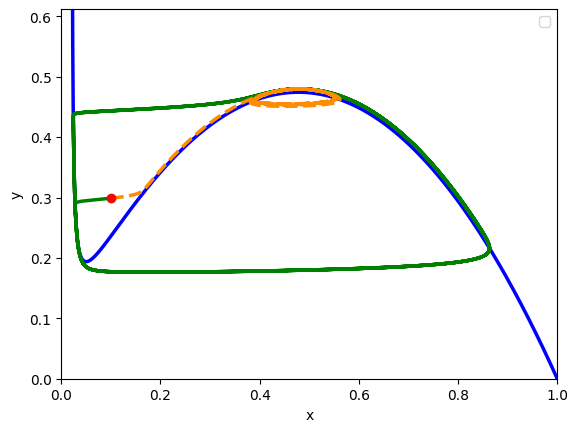} \\
  $(c_1)$\ \  $ f=2.18$ & $(c_6)$\ \  $f=0.572$ \\
  \includegraphics[width=0.45\textwidth]{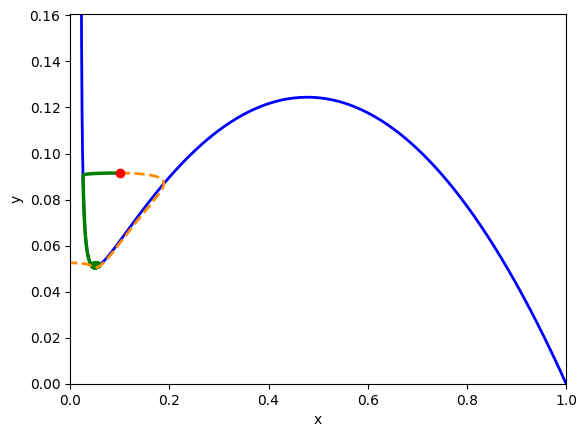} & \includegraphics[width=0.45\textwidth]{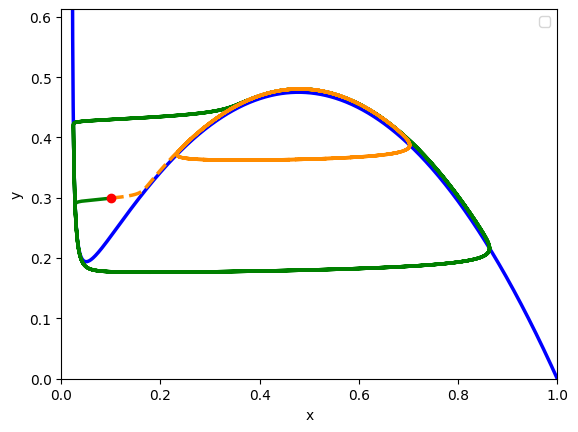} \\
  $(c_2)$\ \  $ f=2.173515959676$ & $(c_7)$\ \  $f=0.57109145$ \\
  \includegraphics[width=0.45\textwidth]{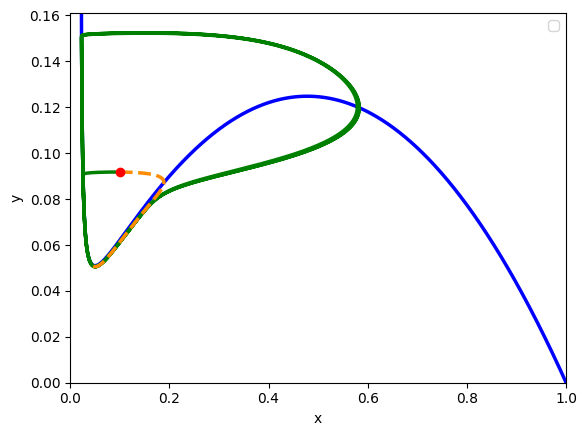} & \includegraphics[width=0.45\textwidth]{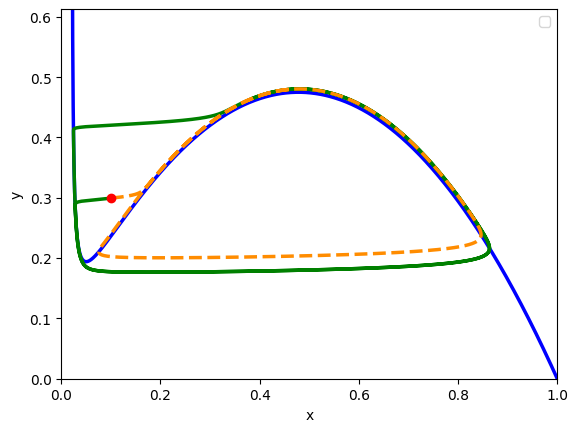} \\
  $(c_3)$\ \  $ f=2.1735159596753$ & $(c_8)$\ \  $f=0.5710914362323$ \\
  \includegraphics[width=0.45\textwidth]{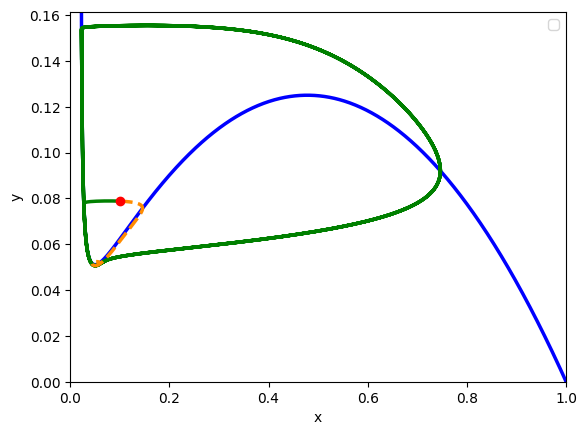} & \includegraphics[width=0.45\textwidth]{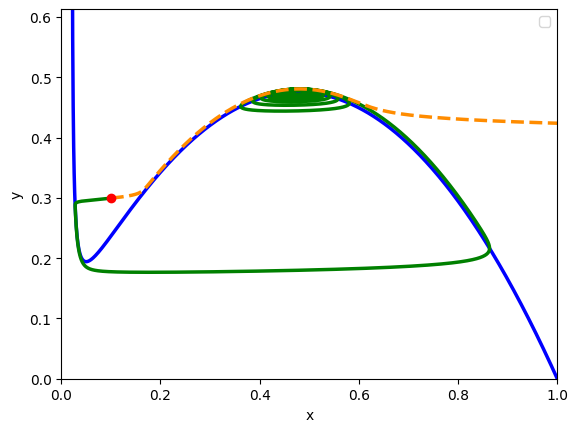} \\
  $(c_4)$\ \  $ f=2.17$ & $(c_9)$\ \  $f=0.571$ \\
  \includegraphics[width=0.45\textwidth]{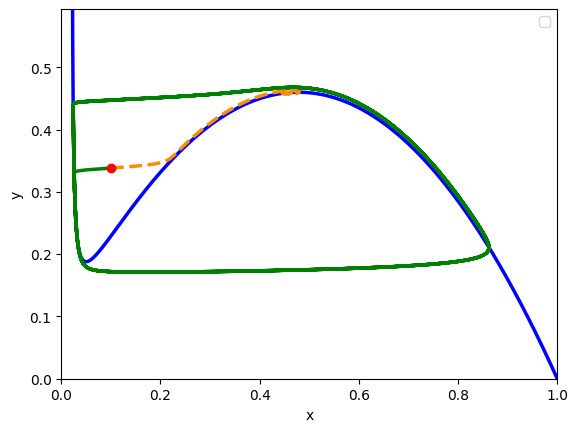} & \includegraphics[width=0.45\textwidth]{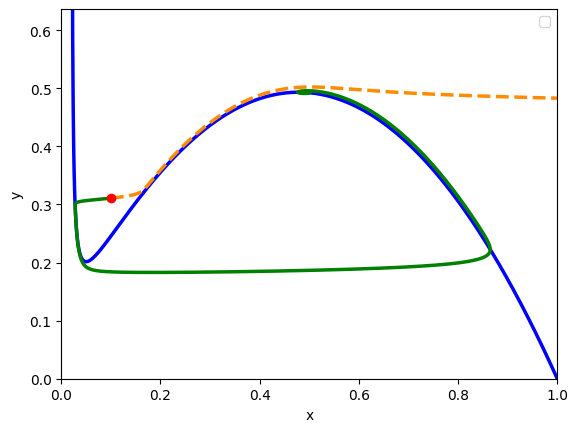} \\
  $(c_5)$\ \  $ f=0.59$ & $(c_{10})$\ \  $f=0.5$ \\
  \caption{Numerical simulation in the case that the sign of $A$ is negative at $E_*=m$ and is positive at $E_*=M$. The green curve represents the positive orbit, and the orange curve represents the negative orbit. $(c_6)-(c_8)$ shows the shrinking of the outer stable limit cycle and the expanding of the inner unstable limit cycle with the decrease of $f$. } \label{sim2} \\
\end{longtable}



\section*{Acknowledgments}

This work is partially supported by National Key R$\&$D Program of China grant
number 2022YFA1005900.

The third author is also partially supported by NNSF of China grant numbers 12071284 and 12161131001.

\end{document}